\documentclass[12pt]{article}
\usepackage[english]{babel}
\usepackage{amssymb}

\newcommand{\R}{{\mathbb R}}
\newcommand{\N}{{\mathbb N}}
\newcommand{\beq}{\begin{equation}}
\newcommand{\ee}{\end{equation}}

\renewcommand{\d}{\partial}

\newtheorem{thm}{Theorem}[section]
\newtheorem{defn}[thm]{Definition}

\newtheorem{example}[thm]{Example}

\newcommand{\om}{\omega}

\newcounter{e}
\newcommand{\reff}[1]{(\ref{#1})}

\title{Fredholm Property of Nonlocal Problems for Integro-Differential Hyperbolic Systems} 

\newcounter{thesame}
\author{
I. Kmit
\ \ \ R. Klyuchnyk\\
{\small
Institute of Mathematics, Humboldt University of Berlin and}
\\
{\small
Institute for Applied Problems of Mechanics and Mathematics,
}
\\
{\small
Ukrainian National Academy of Sciences}
\\
{\small   E-mail:
{\tt kmit@mathematik.hu-berlin.de}}
\\[5mm]
{\small
Institute for Applied Problems of Mechanics and Mathematics,
}
\\
{\small
Ukrainian National Academy of Sciences}
\\
{\small   E-mail:
{\tt roman.klyuchnyk@gmail.com}}
}

\date{}

\begin{document}

\maketitle

\begin{abstract}
The paper concerns nonlocal time-periodic boundary value problems for first-order integro-differential hyperbolic systems with boundary inputs. The systems are subjected to integral boundary conditions.
 Under natural regularity assumptions on the data, it is  shown that the problems display completely non-resonant behavior and satisfy the Fredholm alternative in the spaces of continuous and time-periodic functions.
\end{abstract}



\section{Introduction}\label{sec:intr} 

\renewcommand{\theequation}{{\thesection}.\arabic{equation}}
\setcounter{equation}{0}

\subsection{Motivation}\label{sec:setting}

We consider the following first-order integro-differential hyperbolic system in one space variable
\begin{equation}
 \begin{array}{ll}
 \displaystyle\partial_{t}u_j
 +a_j(x,t)\partial_{x}u_j
 +\sum_{k=1}^{n}b_{jk}(x,t)u_k
 +\sum_{k=1}^{n}\int_{0}^{x}g_{jk}(y,t)u_k(y,t)dy &  \\ [2mm]
\qquad =\displaystyle \sum_{k=m+1}^{n}h_{jk}(x,t)u_k(0,t)
 +\sum_{k=1}^{m}h_{jk}(x,t)u_k(1,t)
 + f_j(x,t),
 \;\;\; x\in(0,1),\; j\le n,&
 \end{array}
 \label{f1}
 \end{equation}
 subjected to periodic conditions in time 
 \begin{equation}\label{f1*}
 u_{j}(x,t)=u_{j}(x,t+2\pi), \;\;\; j\le n, 
 \end{equation}
 and integral boundary conditions in space
 \begin{equation}\label{f2}
 \begin{array}{ll}
   u_{j}(0,t)=
 \displaystyle\sum_{k=1}^{n}\int_{0}^{1}r_{jk}(x,t)u_{k}(x,t)\, dx,
 \;\;\; 1\le j\le m,  \\ [2mm]
  u_{j}(1,t)=
 \displaystyle\sum_{k=1}^{n}\int_{0}^{1}r_{jk}(x,t)u_{k}(x,t)\, dx,
 \;\;\; m<j\le n,
 \end{array}
 \end{equation}
where $0\le m\le n$ are positive integers. 

Note that the boundary terms   
$u_{m+1}(0,t), \dots, u_n(0,t)$ and
$u_1(1,t), \dots,$ $u_m(1,t)$ 
contribute into the differential system \reff{f1}, while the boundary terms   
$u_{1}(0,t), \dots,$ $u_m(0,t)$ and $u_{m+1}(1,t)$, 
$\dots,u_n(1,t)$
contribute into the boundary 
conditions \reff{f2}. 
In this form, which is  motivated by 
applications,  the problem has been studied in \cite{KrsSm,SanNak}.

The  Volterra integral terms in  \reff{f1} are motivated by 
the aforementioned applications (see, e.g., \cite{KrsSm,SanNak}).
As  it will  be seen from our proof of  Theorem \ref{thm:th12}, 
our analysis applies also to the case when  these terms
are replaced by the Fredholm integral terms.

In general, systems of the type \reff{f1}, \reff{f2} model a broad range of physical problems such as traffic 
flows, chemical reactors and heat exchangers \cite{SanNak}. They
are also used to describe problems
of population dynamics (see, e.g., \cite{eft,Keyfitz,mogulru,webb} and references therein) and polymer 
rheology \cite{engl}. Moreover, they  appear in the study of optimal boundary control problems 
\cite{KrsSm,Naka,SanNak,coron}.

Establishing a Fredholm property  is a first step in developing a theory 
of local smooth continuation  \cite{KmRe4} and bifurcation \cite{bueft,cush,KmRe3} for Fredholm hyperbolic operators, in particular, 
such tools as Lyapunov-Schmidt reduction.  
Buono and Eftimie \cite{bueft}  consider  autonomous $2\times 2$ nonlocal hyperbolic systems
in a single space variable, 
describing formation and movement of various animal, cell and bacterial aggregations,
with some biologically motivated integral terms in the 
differential equations. One of the main results in \cite{bueft} is a Fredholm alternative for 
the linearizations at a steady-state, which enables  performing a bifurcation analysis by means of the Lyapunov-Schmidt reduction. 
Here we continue this line of research, establishing the Fredholm property
for a wide range of non-autonomous nonlocal  problems for $(n\times n)$-hyperbolic systems,
with nonlocalities both in the differential equations and in the 
boundary conditions.

We show that the problem \reff{f1}--\reff{f2} demonstrates a completely non-resonant behavior (in other terms, no {\it small divisors} occur). 
More precisely, we prove the Fredholm alternative for 
 \reff{f1}--\reff{f2} under the only assumptions that the coefficients in \reff{f1} and \reff{f2}
are sufficiently smooth and a kind of Levy condition is fulfilled. The proof extends the ideas of 
\cite{KR3,km1} for proving
 the Fredholm alternative  for first-order 
one-dimensional hyperbolic systems with reflection boundary conditions, and also the ideas of 
\cite{Kmit}  for proving a smoothing property for  boundary value hyperbolic problems.
In contrast to \cite{KR3} and \cite{km1}, where  conditions excluding  a resonant behavior
are imposed,
the present Fredholmness result is unconditional, in this respect.

\subsection{Our result}\label{sec:contr}

By $C_{n,2\pi}$ we denote the vector space of all $2\pi$-periodic in $t$ and continuous maps $u:[0,1]\times\R\to\R^{n}$, with the norm
$$
 \|u\|_{\infty}=\max_{j\le n}\max_{x\in [0,1]}\max_{t\in\R}|u_j|.
$$
 Similarly, $C_{n,2\pi}^1$ denotes the Banach space of all $u\in C_{n,2\pi}$ such that $\d_xu,\d_tu\in C_{n,2\pi}$, with the norm 
$$
\|u\|_{1}=\|u\|_{\infty}+\|\partial_x u\|_{\infty}+\|\partial_t u\|_{\infty}.
 $$
For simplicity, we  skip subscript $n$ if $n=1$ and write $C_{2\pi}$ and $C^1_{2\pi}$ for 
$C_{1,2\pi}$ and $C^{1}_{1,2\pi}$, respectively.

We make the following natural assumptions on the coefficients of \reff{f1} and \reff{f2}:
 \begin{eqnarray}
&& \hskip-10mm a_j\in C^1_{2\pi} \mbox{ and } b_{jk}, \d_tb_{jk}, g_{jk}, h_{jk}, r_{jk}, \d_tr_{jk}\in C_{2\pi} \mbox{ for all } j\le n \mbox{ and } k\le n,\label{f4}
\\ [1mm]
&&  \hskip-10mm a_j\neq0  \mbox{ for all } (x,t)\in[0,1]\times\R \mbox{ and } j\le n,\label{f5}
\end{eqnarray}
and 
\beq\label{fz8}
 \begin{array}{ll}
\mbox{for all } 1\le j\neq k\le n \mbox{ there exists } \tilde{b}_{jk}\in C_{2\pi} \mbox{ such that }\\ 
\d_t\tilde b_{jk}  \in C_{2\pi}   \mbox{ and } 
b_{jk}=\tilde{b}_{jk}(a_k-a_j). 
\end{array}
\ee
The assumption \reff{f5} is standard  and means the non-degeneracy of the
hyperbolic system (1.1).  The assumption (1.6) is a kind of the well-known Levy condition 
appearing in various 
aspects of the hyperbolic theory, for instance, for
proving the spectrum-determined growth condition  for semiflows generated 
by initial value problems for hyperbolic systems \cite{Guo,Lichtner,Neves}. 
It plays also a crucial role 
in the Fredholm  analysis of hyperbolic PDEs (see Example \ref{thm:ex1} below).

Given $j\le n$, $x\in[0,1]$, and $t\in\R$, the $j$-th characteristic of \reff{f1} is defined as the solution $\xi\in[0,1]\mapsto\omega_j(\xi,x,t)\in\R$ of the initial value problem
\beq\label{f7}
\d_{\xi}\omega_{j}(\xi,x,t)=\frac{1}{a_j(\xi,\omega_{j}(\xi,x,t))}, \;\;\; \omega_{j}(x,x,t)=t.
\ee
To shorten notation, we will write $\omega_j(\xi)=\omega_j(\xi,x,t)$. In what follows we will use the equalities
\beq \label{f22}
 \partial_x\omega_{j}(\xi)=
 -\frac{1}{a_j(x,t)}\exp{{\int_{\xi}^{x}\left (\frac{\partial_2a_j}{a_j^2}\right )(\eta,\omega_j(\eta))\, d\eta}}
 \ee
and
 \beq \label{f23}
 \partial_t\omega_{j}(\xi)=
 \exp{{\int_{\xi}^{x}\left (\frac{\partial_2a_j}{a_j^2}\right )(\eta,\omega_j(\eta))\, d\eta}},
 \ee
where by $\d_i$ here and below we denote the partial derivative with respect to the $i$-th argument. 
Set 
\beq\label{f8}
\begin{array}{ll}
\displaystyle c_j(\xi,x,t)=
\exp{{\int_{x}^{\xi}\left (\frac{b_{jj}}{a_j}\right )(\eta,\omega_j(\eta))\, d\eta}}, \;\;\;\\ [5mm]
\displaystyle d_j(\xi,x,t)=
\frac{c_j(\xi,x,t)}{a_j(\xi,\omega_j(\xi))},
\end{array}
\ee
and

$$
 x_j=\left\{
 \begin{array}{rl}
 0 &\mbox{if}\ 1\le j\le m,\\
 1 &\mbox{if}\ m<j\le n.
\end{array}
\right.
$$

Integration along the characteristic curves brings the system \reff{f1}--\reff{f2} to the integral form
\begin{eqnarray}\label{f10}
 \displaystyle 
 u_j(x,t)
 &=&\displaystyle c_j(x_j,x,t)\sum_{k=1}^{n}\int_{0}^{1}r_{jk}(\eta,\omega_j(x_j))u_{k}(\eta,\omega_{j}(x_j))\, d\eta   \nonumber  \\
 &&-\displaystyle\sum_{k\neq j}\int_{x_j}^{x}d_{j}(\xi,x,t)b_{jk}(\xi,\omega_j(\xi))u_{k}(\xi,\omega_j(\xi))\, d\xi     \nonumber  \\ 
 &&-\displaystyle\sum_{k=1}^{n}\int_{x_j}^{x}d_{j}(\xi,x,t)\int_{0}^{\xi}g_{jk}(y,\omega_j(\xi))u_{k}(y,\omega_j(\xi))\, dyd\xi \\
 &&+\displaystyle \sum_{k=1}^{n}\int_{x_j}^{x}d_{j}(\xi,x,t)h_{jk}(\xi,\omega_j(\xi))u_{k}(1-x_k,\omega_j(\xi))\, d\xi
 \nonumber \\
 &&+\displaystyle\int_{x_j}^{x}d_{j}(\xi,x,t)f_{j}(\xi,\omega_j(\xi))\, d\xi, \quad j\le n.\nonumber 
\end{eqnarray}
Indeed, let $u$ be a  $C^1$-solution to  \reff{f1}--\reff{f2}. Then,
using  \reff{f1} and  \reff{f7}, for all $j\le n$ we  have
\begin{eqnarray*}
\lefteqn{
\frac{d}{d\xi} u_j(\xi,\om_j(\xi))=
\d_1u_j(\xi,\om_j(\xi))+\frac{\d_2u_j(\xi,\om_j(\xi))}{a_j(\xi,\om_j(\xi))}}\\
&&=\frac{1}{a_j(\xi,\om_j(\xi))}\Biggl(-\sum_{k=1}^n 
b_{jk}(\xi,\om_j(\xi))u_k(\xi,\om_j(\xi))+\displaystyle \sum_{k=m+1}^{n}h_{jk}(\xi,\omega_j(\xi))u_k(0,\omega_j(\xi))
  \\ [2mm]
&&\displaystyle+\sum_{k=1}^{m}h_{jk}(\xi,\omega_j(\xi))u_k(1,\omega_j(\xi))-
\sum_{k=1}^{n}\int_{0}^{\xi}g_{jk}(y,\omega_j(\xi))u_k(y,\omega_j(\xi))\,dy 
+f_j(\xi,\om_j(\xi))\Biggr).
\end{eqnarray*}
This is a linear inhomogeneous ordinary differential equation for the function $u_j(\cdot,\om_j(\cdot,x,t))$, 
and the variation of constants formula (with initial condition at $x_j$) gives
\begin{eqnarray*}
\lefteqn{
u_j(x,t)=u_j(x_j,\om_j(x_j))\exp \int^{x_j}_{x} \left(\frac{b_{jj}}{a_j}\right) (\xi,\om_j(\xi))\,d\xi
-\int^{x_j}_{x}\exp \int_\xi^x\left(\frac{b_{jj}}{a_j}\right)(\eta,\om_j(\eta))\, d \eta} \\
&&\times 
\frac{1}{a_j(\xi,\om_j(\xi))}\Biggl(-\sum_{k\ne j} 
b_{jk}(\xi,\om_j(\xi))u_k(\xi,\om_j(\xi))+\displaystyle \sum_{k=m+1}^{n}h_{jk}(\xi,\omega_j(\xi))u_k(0,\omega_j(\xi))
  \\ [2mm]
&&\displaystyle+\sum_{k=1}^{m}h_{jk}(\xi,\omega_j(\xi))u_k(1,\omega_j(\xi))-
\sum_{k=1}^{n}\int_{0}^{\xi}g_{jk}(y,\omega_j(\xi))u_k(y,\omega_j(\xi))\,dy 
+f_j(\xi,\om_j(\xi))\Biggr)\,
d\xi.
\end{eqnarray*}
Inserting the boundary conditions \reff{f2} and using the notation
\reff{f8},  
we get \reff{f10}, as desired.

\begin{defn}\rm
A function $u\in C_{n,2\pi}$ is called a {\it continuous solution} to \reff{f1}--\reff{f2} if it satisfies \reff{f10}.
\end{defn}
Our result states that either the space of nontrivial solutions to \reff{f1}--\reff{f2} with $f=(f_{1},...,f_{n})=0$ is not empty and has finite dimension or the system \reff{f1}--\reff{f2} has a unique solution for any $f$.

\begin{thm}\label{thm:th12}
Suppose that the conditions \reff{f4}--\reff{fz8} are fulfilled. 
Let $\mathcal{K}$ denote the vector space of all continuous solutions to \reff{f1}--\reff{f2} with $f\equiv0$.
Then

$(i)$ $\dim \mathcal{K}<\infty$ and the vector space of all $f\in C_{n,2\pi}$ such that there exists a continuous solution to \reff{f1}--\reff{f2} is a closed subspace of codimension $\dim \mathcal{K}$ in $C_{n,2\pi}$.

$(ii)$ If $\dim \mathcal{K}=0$, then for any $f\in C_{n,2\pi}$ there exists a unique continuous solution $u$ to \reff{f1}--\reff{f2}.
\end{thm}

\begin{example}\label{thm:ex1}\rm
Consider the following example showing that the condition
 \reff{fz8} plays a crucial role for our result:
\begin{equation}
 \begin{array}{ll}
 \displaystyle\partial_{t}u_1
 +\frac{2}{\pi}\partial_{x}u_1
 -u_2=0 &  \\ [3mm]
\displaystyle\partial_{t}u_2
 +\frac{2}{\pi}\partial_{x}u_2
 +u_1=0,  &
 \end{array}
 \label{f1ex}
 \end{equation} 
 \begin{equation}\label{f3}
 u_{1}(x,t)=u_{1}(x,t+2\pi), \;\;\; u_{2}(x,t)=u_{2}(x,t+2\pi), 
 \end{equation}
 \beq\label{f2ex}
 \begin{array}{ll}
  u_{1}(0,t)=0,  \;\;\;\; 
  u_{2}(1,t)=0. 
  \end{array}
  \ee
This problem is a particular case of \reff{f1}--\reff{f2}
 and satisfies all assumptions of Theorem~\ref{thm:th12} with the exception of \reff{fz8}. 
It is straightforward to check that
$$
u_1=\sin{\frac{\pi}{2}x}\sin{l\left(t-\frac{\pi}{2}x\right)}, \;\;\; 
u_2=\cos{\frac{\pi}{2}x}\sin{l\left(t-\frac{\pi}{2}x\right)}, \quad l\in\N,
$$
are infinitely many linearly independent solutions to the problem \reff{f1ex}--\reff{f2ex} and, therefore, the kernel of the operator of \reff{f1ex}--\reff{f2ex} is infinite dimensional. Thus, the conclusion of Theorem \ref{thm:th12} is not true without \reff{fz8}.

\end{example}

\section{Proof of Theorem \ref{thm:th12}}\label{sec:Fredh}

\renewcommand{\theequation}{{\thesection}.\arabic{equation}}
\setcounter{equation}{0}

Define linear bounded operators $R,B,G,H,F: C_{n,2\pi}\to C_{n,2\pi}$ by
\begin{eqnarray}
 (Ru)_j(x,t)&=&
 c_j(x_j,x,t)\sum_{k=1}^{n}\int_{0}^{1}r_{jk}(\eta,\omega_j(x_j))u_{k}(\eta,\omega_{j}(x_j))\,d\eta,\quad j\le n, \nonumber
\\
 \label{f34}
 (Bu)_j(x,t)&=&
 -\sum_{k\neq j}\int_{x_j}^{x}d_j(\xi,x,t)b_{jk}(\xi,\omega_j(\xi))u_k(\xi,\omega_j(\xi))\, d\xi,
\quad j\le n, \\
\label{f3014}
(Gu)_j(x,t)&=& 
-\sum_{k=1}^{n}\int_{x_j}^{x}\int_{0}^{\xi}d_{j}(\xi,x,t)g_{jk}(y,\omega_j(\xi))
u_{k}(y,\omega_j(\xi))\, dyd\xi,
\quad j\le n, \\
\label{f3024}
 (Hu)_j(x,t)&=&
 \sum_{k=1}^{n}\int_{x_j}^{x}d_{j}(\xi,x,t)h_{jk}(\xi,\omega_j(\xi))u_{k}(1-x_k,\omega_j(\xi))\, d\xi,\quad j\le n,
 \end{eqnarray}
and 
$$
 (Ff)_j(x,t)=
 \int_{x_j}^{x}d_j(\xi,x,t)f_j(\xi,\omega_j(\xi))\, d\xi,\quad j\le n.
$$
 Then the system \reff{f10} can be written as the operator equation
$$
 u=Ru+Bu+Gu+Hu+Fu.
$$

Note that Theorem \ref{thm:th12} says exactly that the operator $I-R-B-G-H: C_{n,2\pi}\to C_{n,2\pi}$ 
is Fredholm of index zero. Nikolsky's criterion  \cite[Theorem XIII.5.2]{KA} says that an operator $I+K$ on a Banach space is Fredholm of index zero whenever $K^2$ is compact. It is interesting to note that the compactness of $K^2$ 
and the identity $I-K^2=(I+K)(I-K)$ imply that the operator $I-K$ is a parametrix of the operator $I+K$ (see \cite{zeid}).

We, therefore, have to show that the operator $K^2: 
C_{n,2\pi}\to C_{n,2\pi}$ for $K^2=(R+B+G+H)^2$ is compact. 
Since the operators $R,B,G$, and $H$ are bounded and the composition of a bounded and a compact operator is 
compact, it is enough to  show that
\beq \label{hrrbbrbg}
\textrm{the operators } \; H, G, R^2, RB, B^2, BR: C_{n,2\pi}\to C_{n,2\pi} \; \textrm{ are compact.}
\ee

We start with the compactness of $H$. By $C_{2\pi}(\R)$ we denote the space of all continuous and $2\pi$-time-periodic maps $v:\R\to\R$. Fix arbitrary $j\le n$ and $k\le n$ and define the operator $H_{jk}\in\mathcal{L}(C_{2\pi}(\R),C_{2\pi})$ by
\begin{equation}\label{hj}
(H_{jk}v)(x,t)
=\int_{x_j}^{x}d_{j}(\xi,x,t)h_{jk}(\xi,\omega_j(\xi))v(\omega_j(\xi))\, d\xi.
\end{equation}
It suffices to show
 the compactness of $H_{jk}$. Change the variable $\xi$ to $z=\omega_j(\xi)$ and denote the inverse 
map by $\xi=\tilde{\omega}_{j}(z)=\tilde{\omega}_{j}(z,x,t)$. Afterwards \reff{hj} reads 
\begin{equation}\label{frfr}
(H_{jk}v)(x,t)
=\int_{\omega_j(x_j)}^{t}d_{j}(\tilde{\omega}_{j}(z),x,t)h_{jk}(\tilde{\omega}_{j}(z),z)a_{j}(\tilde{\omega}_{j}(z),z)v(z)\, dz.
\end{equation}
By the regularity assumption \reff{f4}, the
functions $\omega_j(x_j)$, $\tilde{\omega}_{j}(z)$, $d_{j}(\xi,x,t)$, 
$h_{jk}(x,z)$, and $a_{j}(x,z)$
are continuous in all their arguments and $2\pi$-periodic in $t$ and, hence, are uniformly continuous in $x$ and $t$. 
Then the equicontinuity property of $(H_{jk}v)(x,t)$ for $v$ over a bounded subset of
 $C_{2\pi}(\R)$ straightforwardly follows. Using the Arzela-Ascoli precompactness criterion,
we conclude that $H_{jk}$ and, hence,  $H$ are compact.

Now we consider the operator $G$. Changing the variable $\xi$ to $z=\omega_j(\xi,x,t)$ in \reff{f3014}, we get
\begin{equation}\label{frfrfr}
(Gu)_j(x,t)= 
-\sum_{k=1}^{n}\int_{\omega_j(x_j)}^{t}\int_{0}^{\tilde{\omega}_j(z)}d_{j}(\tilde{\omega}_j(z),x,t)g_{jk}(y,z)a_j(\tilde{\omega}_{j}(z),z)u_{k}(y,z)\, dydz.
\end{equation}
Similarly to the above, the functions $\omega_j(x_j), \tilde{\omega}_j(z), d_j(\tilde{\omega}_j(z),x,t)$, and $a_j(\tilde{\omega}_j(z),z)$ are $2\pi$-periodic in $t$ and uniformly continuous in $x$ and $t$. This entails the equicontinuity property for $(Gu)_j(x,t)$ for $u$ over a bounded subset of $C_{n,2\pi}$. The compactness of $G$ again follows from the Arzela-Ascoli theorem.

We further proceed with the compactness of $R^2$. For $j\le n$ and $k\le n$ define operators $R_{jk}\in\mathcal{L}(C_{2\pi})$ by
$$
(R_{jk}w)(x,t)=c_j(x_j,x,t)\int_{0}^{1}r_{jk}(\eta,\omega_j(x_j))w(\eta,\omega_j(x_j))\,d\eta.
$$
Fix arbitrary $j\le n$, $k\le n$, and $i\le n$. We prove the compactness of the operator 
$R_{jk}R_{ki}$;  the compactness of all other operators contributing into the $R^2$ will follow 
from the same argument. Introduce operators $P_j,Q_{jk} : C_{2\pi}\to C_{2\pi}$ by
\begin{eqnarray}
(P_jw)(x,t)
&=&c_{j}(x_j,x,t)\int_{0}^{1}w(\eta,t)\,d\eta, \\
(Q_{jk}w)(x,t)
&=&r_{jk}(x,\omega_j(x_j))w(x,\omega_j(x_j)).
\end{eqnarray}
Then we have 
$$
R_{jk}=P_jQ_{jk}, \;\;\; R_{ki}=P_kQ_{ki}
$$
and, hence
$$
R_{jk}R_{ki}=P_jQ_{jk}P_kQ_{ki}.
$$
We aim at showing the compactness of $P_jQ_{jk}P_k$, as this and the boundedness of $Q_{ki}$ will entail the compactness of $R_{jk}R_{ki}$. The operator $P_jQ_{jk}P_k$ reads
\begin{equation}\label{ghng}
\begin{array}{lc}
(P_jQ_{jk}P_{k}w)(x,t) = c_j(x_j,x,t)  \\ [2mm]
\displaystyle \quad\times\int_{0}^{1} r_{jk}(\xi,\omega_j(x_j,\xi,t))c_k(x_k,\xi,\omega_j(x_j,\xi,t))\int_{0}^{1}w(\eta,\omega_k(x_k,\xi,t))\, d\eta d\xi.  &
\end{array}
\end{equation}
Changing the variable $\xi$ to $z=\omega_k(x_k,\xi,t)$, we get
\begin{equation}\label{ghng1}
\begin{array}{lr}
(P_jQ_{jk}P_{k}w)(x,t) =  c_j(x_j,x,t)  \\ [2mm]
\displaystyle\quad\times\int_{\omega_k(x_k,0,t)}^{\omega_k(x_k,1,t)}\hskip-5mm r_{jk}(\tilde{\omega}_k(t,x_k,z),z)c_k(x_k,\tilde{\omega}_k(t,x_k,z),z)\int_{0}^{1}\d_{3}\tilde{\omega}_k(t,x_k,z)w(\eta,z)\, d\eta dz,  &
\end{array}
\end{equation}
where
\beq \label{2star}
\partial_3\tilde{\omega}_{k}(\tau,x,t)=
 a_k(x,t)\exp{\int_{\tau}^{t}\partial_1a_k(\tilde{\omega}_k(\rho,x,t),\rho)\, d\rho}.
\ee
Similarly to the above, the compactness of $P_jQ_{jk}P_k$ now immediately follows from the regularity assumption \reff{f4} and the Arzela-Ascoli theorem.

Now we treat the operator 
\begin{eqnarray*}
&(RBu)_{j}(x,t)
=\displaystyle-c_{j}(x_j,x,t)\sum_{k\neq l}\int_{0}^{1}\int_{x_k}^{\eta}r_{jk}(\eta,\omega_j(x_j))
d_{k}(\xi,\eta,\omega_j(x_j)) \nonumber&  \\ [4mm]
&\displaystyle \qquad\times b_{kl}(\xi,\omega_k(\xi,\eta,\omega_j(x_j)))u_{l}(\xi,\omega_k(\xi,\eta,\omega_j(x_j)))\,d\xi d\eta \nonumber&
\end{eqnarray*}
for an arbitrary fixed $j\le n$.
After changing the order of integration we get the equality
\begin{eqnarray*}
&(RBu)_{j}(x,t)
=\displaystyle-c_{j}(x_j,x,t)\sum_{k\neq l}\int_{0}^{1}\int_{\xi}^{1-x_k}r_{jk}(\eta,\omega_j(x_j))d_{k}(\xi,\eta,\omega_j(x_j)) &  \\ [2mm]
&\displaystyle\times b_{kl}(\xi,\omega_k(\xi,\eta,\omega_j(x_j)))u_{l}(\xi,\omega_k(\xi,\eta,\omega_j(x_j)))\,d\eta d\xi.  & 
\end{eqnarray*}
Then we change the variable $\eta$ to $z=\omega_k(\xi,\eta,\omega_j(x_j))$. 
Since the inverse is given by $\eta=\tilde{\omega}_k(\omega_j(x_j),\xi,z)$,  we get
\begin{eqnarray}
\lefteqn{
(RBu)_{j}(x,t)
=}\nonumber\\ &&-\displaystyle c_{j}(x_j,x,t)\sum_{k\neq l}
\displaystyle\int_{0}^{1}\int_{\omega_j(x_j)}^{\omega_k(\xi,1-x_k,\omega_j(x_j))}r_{jk}(\tilde{\omega}_k(\omega_j(x_j),\xi,z),\omega_j(x_j))
\label{rdoc}   \\ [2mm]
&&\times\displaystyle d_{k}(\xi,\tilde{\omega}_k(\omega_j(x_j),\xi,z),\omega_j(x_j))b_{kl}(\xi,z)
\displaystyle \d_{3}\tilde{\omega}_k(\omega_j(x_j),\xi,z)u_{l}(\xi,z)\,dz d\xi,\nonumber
\end{eqnarray}
where $\d_{3}\tilde{\omega}_k(\omega_j(x_j),\xi,z)$ is given by \reff{2star}. The functions $\omega_j(\xi,x,t)$ and the kernels of the integral operators in \reff{rdoc} are continuous and $t$-periodic functions and, hence, are uniformly continuous functions in $x$ and $t$.
This means that we are again in the conditions of the Arzela-Ascoli theorem, as desired.

We proceed to show that $B^2:C_{n,2\pi}\to C_{n,2\pi}$ is compact. By the Arcela-Ascoli theorem, $C_{n,2\pi}^1$ is compactly embedded into $C_{n,2\pi}$. Then the desired compactness property will follow if we show that
\beq \label{f309}
B^2 \;\textrm{maps continuously}\; C_{n,2\pi}\; \textrm{into} \; C_{n,2\pi}^1.
\ee
By using the equalities \reff{f22}, \reff{f23}, and \reff{f34},
the partial derivatives $\partial_xB^2u$, $\partial_{t}B^2u$ exist and are continuous for each $u\in C_{n,2\pi}^1$. Since $C^1_{n,2\pi}$ is dense in $C_{n,2\pi}$, the desired condition \reff{f309} will follow from the bound
\beq \label{bbboun} 
\left\|B^2u\right\|_{1}
= O\left(\|u\|_{\infty}\right)\; \textrm{for all} \; u\in C^1_{n,2\pi}.
\ee 
 To prove \reff{bbboun}, for given $j\le n$ and $u\in C_{n,2\pi}^1$, let us consider the following representation for $(B^2u)_j(x,t)$ obtained after the application of the Fubini's theorem:
\beq \label{f311}
(B^2u)_j(x,t)
=\sum_{k\neq j}\sum_{l\neq k}\int_{x_j}^x\int_{\eta}^{x} d_{jkl}(\xi,\eta,x,t)b_{jk}(\xi,\omega_j(\xi))u_l(\eta,\omega_k(\eta,\xi,\omega_j(\xi)))\, d\xi d\eta,
\ee
where 
\beq \label{f311d}
d_{jkl}(\xi,\eta,x,t)=d_j(\xi,x,t)d_{k}(\eta,\xi,\omega_{j}(\xi))b_{kl}(\eta,\omega_{k}(\eta,\xi,\omega_j(\xi))).
\ee
The estimate $\left\|B^2u\right\|_{\infty}=O\left(\|u\|_{\infty}\right)$  is obvious.
Since 
$$
(\d_t+a_j(x,t)\d_x)\varphi(\omega_j(\xi,x,t))=0
$$
for all $j\le n, \varphi\in C^1(\R), x,\xi\in[0,1]$, and $t\in\R$, one can easily check that
$$
\|[(\d_t+a_j(x,t)\d_x)(B^2u)_j]\|_{\infty} = O\left(\|u\|_{\infty}\right)
\mbox{ for all } j\le n \mbox{ and } u\in C^1_{n,2\pi}.
$$
Hence the estimate 
$\left\|\d_xB^2u\right\|_{\infty}=O\left(\|u\|_{\infty}\right)$ will follow from the following one:
\beq \label{f31jhr1}
\|\d_tB^2u||_{\infty}= O\left(\|u\|_{\infty}\right).
\ee
In order to prove \reff{bbboun}, we are therefore reduced to prove \reff{f31jhr1}. To this end, we start with the following consequence of \reff{f311}:
\begin{eqnarray}
\lefteqn{
\d_t[(B^2u)_j(x,t)]}
\nonumber\\ &&
=\displaystyle\sum_{k\neq j}\sum_{l\neq k}\int_{x_j}^x\int_{\eta}^{x} \frac{d}{dt}\Bigl[ d_{jkl}(\xi,\eta,x,t)b_{jk}(\xi,\omega_j(\xi))\Bigr] u_l(\eta,\omega_k(\eta,\xi,\omega_j(\xi)))\, d\xi d\eta
\nonumber\\ &&
+\displaystyle\sum_{k\neq j}\sum_{l\neq k}\int_{x_j}^x\int_{\eta}^{x} d_{jkl}(\xi,\eta,x,t) b_{jk}(\xi,\omega_j(\xi))
\nonumber\\ &&
\times
\d_t\omega_k(\eta,\xi,\omega_j(\xi))\d_t\omega_j(\xi)\d_2u_l(\eta,\omega_k(\eta,\xi,\omega_j(\xi))) \,d\xi d\eta. \nonumber
\end{eqnarray}
Let us transform the second summand. Using \reff{f7}, \reff{f22}, and \reff{f23}, we get
\begin{eqnarray}
\lefteqn{
\frac{d}{d\xi} u_l(\eta,\omega_k(\eta,\xi,\omega_j(\xi)))} \nonumber \\ &&
=\Bigl[\d_x\omega_k(\eta,\xi,\omega_j(\xi))+\d_t\omega_k(\eta,\xi,\omega_j(\xi))\d_{\xi}\omega_j(\xi)\Bigr] \d_2u_l(\eta,\omega_k(\eta,\xi,\omega_j(\xi)))
\label{eqwn}
\\ &&
=\left ( \frac{1}{a_j(\xi,\omega_j(\xi))}-\frac{1}{a_k(\xi,\omega_j(\xi))}\right ) \d_t\omega_k(\eta,\xi,\omega_j(\xi))\d_2u_l(\eta,\omega_k(\eta,\xi,\omega_j(\xi))) \nonumber.
\end{eqnarray}
Therefore,
\begin{eqnarray}
\lefteqn{ 
 b_{jk}(\xi,\omega_j(\xi))\d_t\omega_k(\eta,\xi,\omega_j(\xi))\d_2u_l(\eta,\omega_k(\eta,\xi,\omega_j(\xi)))}
\nonumber \\ &&
 =\displaystyle a_j(\xi,\omega_j(\xi))a_k(\xi,\omega_j(\xi))\tilde{b}_{jk}(\xi,\omega_j(\xi))\frac{d}{d\xi} u_l(\eta,\omega_k(\eta,\xi,\omega_j(\xi))), \label{f313}
\end{eqnarray}
where the functions $\tilde{b}_{jk}\in C_{2\pi}$ are fixed to satisfy \reff{fz8}. Note that $\tilde{b}_{jk}$ are not uniquely defined by \reff{fz8} for $(x,t)$ with $a_{j}(x,t)=a_{k}(x,t)$. Nevertheless, as it follows from \reff{eqwn}, the right-hand side (and, hence, the left-hand side of \reff{f313}) do not depend on the choice of $\tilde{b}_{jk}$, since $\frac{d}{d\xi}u_{l}(\eta,\omega_{k}(\eta,\xi,\omega_{j}(\xi)))=0$ if $a_{j}(x,t)=a_{k}(x,t)$. 

Write
$$
\tilde{d}_{jkl}(\xi,\eta,x,t)
=d_{jkl}(\xi,\eta,x,t)\d_t\omega_j(\xi)a_k(\xi,\omega_j(\xi))a_{j}(\xi,\omega_j(\xi))\tilde{b}_{jk}(\xi,\omega_j(\xi)),
$$
where $d_{jkl}$ are introduced by \reff{f311d} and \reff{f8}. Using \reff{f7} and \reff{f22}, 
we see that the function $\tilde{d}_{jkl}(\xi,\eta,x,t)$ is $C^1$-regular in $\xi$ due to regularity assumptions \reff{f4} and \reff{fz8}. Similarly,
using \reff{f23}, we see that the functions $d_{jkl}(\xi,\eta,x,t)$ and $b_{jk}(\xi,\omega_j(\xi))$ are $C^1$-smooth in $t$.

By \reff{f313} we have
\begin{eqnarray*}
\lefteqn{
(\d_tB^2u)_j(x,t)}\\  &&
= \displaystyle \sum_{k\neq j}\sum_{l\neq k}\int_{x_j}^x\int_{\eta}^{x} \frac{d}{dt}[ d_{jkl}(\xi,\eta,x,t)b_{jk}(\xi,\omega_j(\xi))] u_l(\eta,\omega_k(\eta,\xi,\omega_j(\xi)))\, d\xi d\eta
\\  &&
+\displaystyle \sum_{k\neq j}\sum_{l\neq k}\int_{x_j}^x\int_{\eta}^{x}\tilde{d}_{jkl}(\xi,\eta,x,t)\frac{d}{d\xi} u_l(\eta,\omega_k(\eta,\xi,\omega_j(\xi)))\, d\xi d\eta
\\  &&
=\displaystyle \sum_{k\neq j}\sum_{l\neq k}\int_{x_j}^x\int_{\eta}^{x} \frac{d}{dt} [d_{jkl}(\xi,\eta,x,t)b_{jk}(\xi,\omega_j(\xi))] u_l(\eta,\omega_k(\eta,\xi,\omega_j(\xi)))\, d\xi d\eta
\\  &&
-\displaystyle \sum_{k\neq j}\sum_{l\neq k}\int_{x_j}^x\int_{\eta}^{x}\d_{\xi}\tilde{d}_{jkl}(\xi,\eta,x,t)u_l(\eta,\omega_k(\eta,\xi,\omega_j(\xi)))\, d\xi d\eta
\\  &&
+\displaystyle \sum_{k\neq j}\sum_{l\neq k}\int_{x_j}^x\left [\tilde{d}_{jkl}(\xi,\eta,x,t) u_l(\eta,\omega_k(\eta,\xi,\omega_j(\xi)))\right ]_{\xi=\eta}^{\xi=x}\, d\eta.
\end{eqnarray*}
The desired estimate \reff{f31jhr1} now easily follows from the assumptions \reff{f4}--\reff{fz8}.

Returning back to \reff{hrrbbrbg}, it remains to prove that the operator $BR:C_{n,2\pi}\to C_{n,2\pi}$ is compact. By the definitions of $B$ and $R$,
\begin{equation}\label{dro}
\begin{array}{cc}
(BRu)_{j}(x,t)
=-\displaystyle \sum_{k\neq j}\sum_{l=1}^{n}\int_{0}^{1}\int_{x_j}^{x}d_{j}(\xi,x,t)b_{jk}(\xi,\omega_j(\xi))c_{k}(x_k,\xi,\omega_j(\xi))
&  \\
\displaystyle \times r_{kl}(\eta,\omega_k(x_k,\xi,\omega_j(\xi)))u_{l}(\eta,\omega_k(x_k,\xi,\omega_j(\xi)))\,d\xi d\eta, \;\;\; j\le n.    & 
\end{array}
\end{equation}
The integral operators in \reff{dro} are similar to those in \reff{f311} and, therefore,
the proof of the compactness of $BR$ follows along the same line as the proof of the compactness of $B^2$. The proof of Theorem \ref{thm:th12} is complete.

\section*{Acknowledgments}

The second author was supported by the BMU-MID Erasmus Mundus Action 2 grant. He expresses his 
gratitude to the Applied Analysis group at the Humboldt University of Berlin for its kind hospitality.

\end{document}